\documentclass[12pt,oneside, italian]{article}
\usepackage{latexsym}
\usepackage{color}
\usepackage{amssymb,amsthm,amsmath}
\usepackage{graphicx}
\usepackage{multirow}

\newtheorem{definition}{Definition}

\begin{document}
\begin{center}
{\bf \LARGE  Classification of minimal\\ 
\vspace*{-0.3 cm}
$1$-saturating sets in $PG(2,q)$, $q\leq 23$\\}
\vspace*{0.4 cm}
{\large D. Bartoli, S. Marcugini and F. Pambianco}\\
{\small Dipartimento di Matematica e Informatica,\\
Universit\`{a} degli Studi di Perugia, \\
Via Vanvitelli 1, 06123 Perugia Italy \\
e-mail: \{daniele.bartoli, gino, fernanda\}@dmi.unipg.it}
\end{center}
\setlength{\baselineskip}{15.5 pt}
{\bf Abstract} 
{\small Minimal $1-$saturating sets in the projective plane $PG(2,q)$ are
considered. They correspond to covering codes which can be applied to many branches of combinatorics and information theory, as data compression, compression with distortion, broadcasting in interconnection network, write-once memory or steganography (see \cite{Coh} and \cite{BF2008}). The full classification of all the minimal $1$-saturating sets in $%
PG(2,9)$ and $PG(2,11)$ and the classification of minimal $1$-saturating sets of smallest size in $PG(2,q)$, $16\leq q\leq 23$ are given. These results have been found using a computer-based exhaustive search that exploits projective equivalence properties.}
\vspace*{0.3 cm}\\
{\bf \large 1 Introduction}\\
\vspace*{-0.3 cm}\\
Let $PG(n,q)$ be the $n$-dimensional projective space over the Galois field $%
GF(q)$. For an introduction to such spaces and the geometrical objects
therein, see \cite{Hirs} - \cite{HirsThas}.
\begin{definition}
A point set $S$ in the space $PG(n,q)$ is $\varrho $-saturating if $\varrho $
is the least integer such that for any point $x\in PG(n,q)$ there exist $%
\varrho +1$ points in $S$ generating a subspace of $PG(n,q)$ in which $x$
lies.
\end{definition}
\begin{definition}
\cite{Ughi} A $\varrho $-saturating set of $l$ points is called minimal if
it does not contain a $\varrho $-saturating set of $l-1$ points.
\end{definition}
A $q$-ary linear code with codimension $r$ has covering radius $R$ if every $%
r$-positional $q$-ary column is equal to a linear combination of $R$ columns
of a parity check matrix of this code and $R$ is the smallest value with
such property. For an introduction to coverings of vector spaces over finite
fields and to the concept of code covering radius, see \cite{Coh}.\\
The points of a $\varrho -$saturating set in $PG(n,q)$ can be considered as
columns of a parity check matrix of a $q$-ary linear code with codimension $%
n+1$. So, in terms of the coding theory, a $\varrho -$saturating $l$-set in $%
PG(n,q)$ corresponds to a parity check matrix of a $q$-ary linear code with
length $l$, codimension $n+1$, and covering radius $\varrho +1$ \cite{Dav95},%
\cite{DavO2},\cite{Janwa}. Such code is denoted by an $[l,\,l-(n+1)]_{q}(%
\varrho +1)$ code. Covering codes can be applied to many branches of combinatorics and information theory, as data compression, compression with distortion, broadcasting in interconnection network, write-once memory or steganography (see \cite{Coh} and \cite{BF2008}).\\
Note that a $\varrho $-saturating set in $PG(n,q),\,\varrho +1\le n$, can
generate an infinite family of $\varrho $-saturating sets in $PG(N,q)$ with $%
N=n+(\varrho +1)m,\,m=1,2,3,\ldots $, see \cite[Chapter 5.4]{Coh},\cite
{Dav95},\cite[Example 6]{DavNBCR2} and references therein, where such
infinite families are considered as linear codes with covering radius $%
\varrho +1$.\\
In the projective plane $PG(2,q)$ over the Galois field $GF(q)$, an $n-$arc
is a set of $n$ points no $3$ of which are collinear. An $n$-arc is called
complete if it is not contained in an $(n+1)$-arc of the same projective
plane. The complete arcs of $PG(2,q)$ are examples of minimal $1$-saturating
sets, but there are minimal $1$-saturating sets that are not complete arcs.
Properties of the $\varrho $-saturating sets in $PG(n,q)$ are presented in \cite{1sat}.\\
\vspace*{-0.2 cm}\\
{\bf \large 2 The computer search for the non-equivalent minimal $1-$saturating
sets}\\
\vspace*{-0.3 cm}\\
Our goal is to determine the classification of saturating sets up to projective equivalence in $PG(2,q)$. The problem of finding non-equivalent geometrical structures is very popular in literature (see \cite{Hirs}, \cite{Hirs2}, \cite{HirsSt}, \cite{HirsThas}). In \cite{MarPam03} the full classification of minimal $1$-saturating sets in $PG(2,q)$, $q\leq 8$, the classification of minimal $1$-saturating sets in $PG(2,q)$ of smallest size for $9\leq q\leq 13$ and the determination of the smallest size of minimal $1$-saturating sets in $PG(2,16)$ are presented.\\
In this work we perform an exhaustive search using a backtracking algorithm. Some strategies have to be used to reduce the search space, as in this kind of problems there are many equivalent parts of the search space and a large number of copies of equivalent solutions could be found. The program starts classifying the sets in $PG(2,q)$ containing the projective frame, until a certain size $k$. We only searched for  minimal $1$-saturating sets containing a projective frame, since the following theorem holds.\\

\vspace*{-0.5 cm}
{\bf Theorem 1} \emph{In $PG(2,q)$ there exists a unique minimal $1$-saturating set not containing a projective frame.  It consists of a whole line and an external point. Its stabilizer has size $\frac{|PGL(3,q)|}{q^2(q^2+q+1)}$ $\left(or \text{ }\frac{|P\Gamma L(3,q)|}{q^2(q^2+q+1)}\right)$}.\\
\vspace*{-0.5 cm}

Then the sets of size $k$ are extended using backtracking. During the backtracking some information obtained during the classification phase is used to further prune the search space. The sets are tested for the saturating property and the minimality condition. See \cite{PhdTesi} for a detailed description.\\
The following tables present the results obtained. In particular we perform a full classification of minimal $1$-saturating sets in $PG(2,9)$ and $PG(2,11)$ and the classification of minimal $1$-saturating sets of smallest size in $PG(2,q)$ with $16\leq q \leq 23$.\\ 
We found no examples of minimal $q+2$-saturating sets in $PG(2,9)$ and $PG(2,11)$ containing the projective frame and then the unique example is that one described above. In the following table we describe the obtained results, in particular the type of the stabilizer of the minimal $1$-saturating sets of size $k$. With the symbol $G_{i}$ we denote a group of order $i$; for the other symbols we refer to \cite{librogruppi}. When complete arcs exist, their number is indicated in bold font. 

{\scriptsize
\vspace*{1 cm}
\begin{center}
\tabcolsep=1 mm
\begin{tabular}{|c|c||  lllll |}
\hline
\multirow{14}{*}{$PG(2,9)$}&\multirow{1}{*}{$k=6$} & {\bf $G_{120}$}: {\bf 1}&&&&\\
&&&&&&\\
&\multirow{1}{*}{$k=7$} & $\mathbb{Z}_{4}$: 1 &  $G_{42}$: {\bf 1}&$G_{120}$: 1&&\\
&&&&&&\\
&\multirow{2}{*}{$k=8$} & $\mathbb{Z}_{1}$: 88&$\mathbb{Z}_{2}$: 52&$\mathbb{Z}_{2}\times \mathbb{Z}_{2}$: 11&$\mathcal{S}_{3}$: 1&$\mathbb{Z}_{2}\times \mathbb{Z}_{4}$: 1\\
&& $\mathcal{D}_{4}$: 1&$\mathcal{D}_{6}$: 3&$G_{16}$: 1+{\bf 1}& $G_{24}$: 2& $G_{48}$: 1\\
&&&&&&\\
&\multirow{2}{*}{$k=9$} &  $\mathbb{Z}_{1}$: 667&$\mathbb{Z}_{2}$: 87&$\mathbb{Z}_{3}$:9&$\mathbb{Z}_{2}\times \mathbb{Z}_{2}$: 4& $\mathcal{S}_{3}$: 2\\
&&$\mathcal{D}_{4}$: 1&$\mathcal{D}_{6}$: 1&$G_{16}$: 1&  $G_{48}$: 1&\\
&&&&&&\\
&\multirow{2}{*}{$k=10$} & $\mathbb{Z}_{1}$: 58&$\mathbb{Z}_{2}$: 22&$\mathbb{Z}_{4}$:5&$\mathbb{Z}_{2}\times \mathbb{Z}_{2}$: 4&$\mathcal{D}_{4}$: 2\\
&&$G_{16}$: 1&  $G_{20}$: 1&  $G_{32}$: 1+{\bf 1}&  $G_{1440}$: 1&\\
&&&&&&\\
&\multirow{1}{*}{$k=11$} &  $G_{11520}$: 1&&&&\\
\hline
\hline
\multirow{18}{*}{$PG(2,11)$}&\multirow{1}{*}{$k=7$} & {\bf $\mathbb{Z}_{7}\rtimes\mathbb{Z}_{3}$}: {\bf 1}&&&&\\
&&&&&&\\
&\multirow{2}{*}{$k=8$} & $\mathbb{Z}_{1}$: 22 &  $\mathbb{Z}_{2}$: 26+{\bf 5}&$\mathbb{Z}_{2}\times\mathbb{Z}_{2}$: 2+{\bf 1 }&$\mathcal{D}_{4}$: 1 + {\bf 1 }&$\mathcal{D}_{5}$: {\bf 1}\\
&&$G_{16}$: { \bf 1}&&&&\\
&&&&&&\\
&\multirow{2}{*}{$k=9$} & $\mathbb{Z}_{1}$: 10686&$\mathbb{Z}_{2}$: 265+{\bf 1}&$\mathbb{Z}_{3}$: 40 +{\bf 1}&
$\mathbb{Z}_{4}$: 2 & $\mathbb{Z}_{2}\times \mathbb{Z}_{2}$: 3 \\
&&$\mathcal{S}_{3}$: 10+{\bf 1}& $\mathbb{Z}_{10}$: 1 & $\mathcal{Q}_6$ : 1&&\\
&&&&&&\\
&\multirow{3}{*}{$k=10$} & $\mathbb{Z}_{1}$: 115731&$\mathbb{Z}_{2}$: 1332&$\mathbb{Z}_{3}$: 31&$\mathbb{Z}_{4}$: 15&$\mathbb{Z}_{2}\times \mathbb{Z}_{2}$: 13\\
&&$\mathbb{Z}_{5}$: 2 &$\mathcal{S}_{3}$: 8&$\mathcal{D}_{4}$: 2&$\mathcal{D}_{5}$: 2&$\mathbb{Z}_{10}$: 1 \\
&&$\mathcal{Q}_6$ :1 &$G_{60}$:{\bf 1} &&&\\
&&&&&&\\
&\multirow{1}{*}{$k=11$} & $\mathbb{Z}_{1}$: 30802&$\mathbb{Z}_{2}$: 147&
$\mathbb{Z}_{4}$: 1 & $\mathbb{Z}_{2}\times \mathbb{Z}_{2}$: 3 &$\mathcal{D}_{4}$: 3\\
&&&&&&\\
&\multirow{2}{*}{$k=12$} & $\mathbb{Z}_{1}$: 119&$\mathbb{Z}_{2}$: 7&$\mathbb{Z}_{3}$: 5&
$\mathcal{S}_{3}$: 1 & $\mathcal{Q}_{6}$: 1 \\
&&$G_{20}$: 1& $G_{1320}$: {\bf 1} &&&\\
&&&&&&\\
&\multirow{1}{*}{$k=13$} & $G_{13200}$: 1&&&&\\
\hline
\hline
\multirow{6}{*}{$q=16$}&
\multirow{1}{*}{$k=9$} & $\mathbb{Z}_{3}$: {\bf 1} &  $\mathbb{Z}_{6}$: {\bf 1}&$\mathcal{D}_{6}$: 1&$\mathcal{G}_{54}$: 1&\\
&&&&&&\\
&\multirow{4}{*}{$k=10$} & $\mathbb{Z}_{1}$: 7744+{\bf 342}&$\mathbb{Z}_{2}$: 699+{\bf 130}&$\mathbb{Z}_{3}$: 3&$\mathbb{Z}_{4}$: 12+{\bf 8}&$\mathbb{Z}_{2}\times  \mathbb{Z}_{2}$: 27+{\bf 4}\\
&&$\mathbb{Z}_{6}$: 2&$\mathcal{S}_{3}$: 4+{\bf 3}&$\mathcal{D}_{4}$: 8 & $\mathbb{Z}_{2}\times \mathbb{Z}_{2}\times  \mathbb{Z}_{2}$: 18+{\bf 10} &$\mathcal{Q}_{6}$: 1\\
&&$\mathbb{Z}_{4}\times \mathbb{Z}_{4}$: 1&$G_{16}$: 4+{\bf 3}&$G_{20}$: {\bf 1}&$G_{24}$: 1&$G_{32}$: 1\\\
&&$G_{48}$: 1&&&&\\
\hline
\hline
\multirow{3}{*}{$q=17$}&
\multirow{3}{*}{$k=10$} & $\mathbb{Z}_{1}$: 2591+{\bf 341}&$\mathbb{Z}_{2}$: 460+{\bf 179}&$\mathbb{Z}_{3}$: 8+{\bf 10}&$\mathbb{Z}_{4}$: 4+{\bf 7}&$\mathbb{Z}_{2}\times  \mathbb{Z}_{2}$: 5+{\bf 8}\\
&&$\mathcal{S}_{3}$: 7+{\bf 9}&$\mathcal{D}_{4}$: 4 &$\mathcal{Q}_{4}$: {\bf 1}&$\mathcal{Q}_{6}$: {\bf 2}& $G_{16}$: 1+{\bf 1}\\
&&$G_{18}$: {\bf 1}&$G_{24}$: {\bf 1}&&&\\
\hline
\hline
\multirow{2}{*}{$q=19$}&
\multirow{2}{*}{$k=10$} & $\mathbb{Z}_{1}$: 1+{\bf 1}&$\mathbb{Z}_{2}$: 6+{\bf 18}&$\mathbb{Z}_{3}$: {\bf 1}&$\mathbb{Z}_{4}$: {\bf 1}&$\mathbb{Z}_{2}\times  \mathbb{Z}_{2}$: {\bf 2}\\
&&$\mathcal{S}_{3}$: {\bf 2}&$\mathcal{D}_{5}$: {\bf 2}&$\mathcal{Q}_{6}$: {\bf 1}&$G_{60}$: {\bf 1}&\\
\hline
\hline
\multirow{1}{*}{$q=23$}&
\multirow{1}{*}{$k=10$} & $\mathcal{S}_{3}$: {\bf 1}&&&&\\
\hline
\end{tabular}
\end{center}
}

\end{document}